%
%
%
%


\magnification=1200
\pretolerance=500 \tolerance=1000 \brokenpenalty=5000
\hsize=12.5cm   
\vsize=19cm
\hoffset=0.4cm
\voffset=1cm   
\parskip3pt plus 1pt
\parindent=0.6cm
\let\sl=\it
\def\\{\hfil\break}


\font\seventeenbf=cmbx10 at 17.28pt

\font\twelvebf=cmbx10 at 12pt
\font\eightbf=cmbx8
\font\sixbf=cmbx6

\font\eighti=cmmi8
\font\sixi=cmmi6

\font\eightrm=cmr8
\font\sixrm=cmr6

\font\eightsy=cmsy8
\font\sixsy=cmsy6

\font\eightit=cmti8
\font\eighttt=cmtt8
\font\eightsl=cmsl8

\font\seventeenbsy=cmbsy10 at 17.28pt

\font\twelvebsy=cmbsy10 at 12pt
\font\tenbsy=cmbsy10
\font\eightbsy=cmbsy8
\font\sevenbsy=cmbsy7
\font\sixbsy=cmbsy6
\font\fivebsy=cmbsy5

\font\tenmsa=msam10

\font\sevenmsa=msam7
\font\fivemsa=msam5
\newfam\msafam
  \textfont\msafam=\tenmsa
  \scriptfont\msafam=\sevenmsa
  \scriptscriptfont\msafam=\fivemsa

\font\tenmsb=msbm10
\font\eightmsb=msbm8
\font\sevenmsb=msbm7
\font\fivemsb=msbm5
\newfam\msbfam
  \textfont\msbfam=\tenmsb
  \scriptfont\msbfam=\sevenmsb
  \scriptscriptfont\msbfam=\fivemsb
\def\Bbb{\fam\msbfam\tenmsb}

\font\tenCal=eusm10
\font\sevenCal=eusm7
\font\fiveCal=eusm5
\newfam\Calfam
  \textfont\Calfam=\tenCal
  \scriptfont\Calfam=\sevenCal
  \scriptscriptfont\Calfam=\fiveCal
\def\Cal{\fam\Calfam\tenCal}

\font\teneuf=eusm10
\font\teneuf=eufm10
\font\seveneuf=eufm7
\font\fiveeuf=eufm5
\newfam\euffam
  \textfont\euffam=\teneuf
  \scriptfont\euffam=\seveneuf
  \scriptscriptfont\euffam=\fiveeuf
\def\euf{\fam\euffam\teneuf}

\font\seventeenbfit=cmmib10 at 17.28pt

\font\twelvebfit=cmmib10 at 12pt
\font\tenbfit=cmmib10
\font\eightbfit=cmmib8
\font\sevenbfit=cmmib7
\font\sixbfit=cmmib6
\font\fivebfit=cmmib5
\newfam\bfitfam
  \textfont\bfitfam=\tenbfit
  \scriptfont\bfitfam=\sevenbfit
  \scriptscriptfont\bfitfam=\fivebfit


\catcode`\@=11
\def\eightpoint{%
  \textfont0=\eightrm \scriptfont0=\sixrm \scriptscriptfont0=\fiverm
  \def\rm{\fam\z@\eightrm}%
  \textfont1=\eighti \scriptfont1=\sixi \scriptscriptfont1=\fivei
  \def\oldstyle{\fam\@ne\eighti}%
  \textfont2=\eightsy \scriptfont2=\sixsy \scriptscriptfont2=\fivesy
  \textfont\itfam=\eightit
  \def\it{\fam\itfam\eightit}%
  \textfont\slfam=\eightsl
  \def\sl{\fam\slfam\eightsl}%
  \textfont\bffam=\eightbf \scriptfont\bffam=\sixbf
  \scriptscriptfont\bffam=\fivebf
  \def\bf{\fam\bffam\eightbf}%
  \textfont\ttfam=\eighttt
  \def\tt{\fam\ttfam\eighttt}%
  \textfont\msbfam=\eightmsb
  \def\Bbb{\fam\msbfam\eightmsb}%
  \abovedisplayskip=9pt plus 2pt minus 6pt
  \abovedisplayshortskip=0pt plus 2pt
  \belowdisplayskip=9pt plus 2pt minus 6pt
  \belowdisplayshortskip=5pt plus 2pt minus 3pt
  \smallskipamount=2pt plus 1pt minus 1pt
  \medskipamount=4pt plus 2pt minus 1pt
  \bigskipamount=9pt plus 3pt minus 3pt
  \normalbaselineskip=9pt
  \setbox\strutbox=\hbox{\vrule height7pt depth2pt width0pt}%
  \let\bigf@ntpc=\eightrm \let\smallf@ntpc=\sixrm
  \normalbaselines\rm}
\catcode`\@=12

\def\eightpointbf{%
 \textfont0=\eightbf   \scriptfont0=\sixbf   \scriptscriptfont0=\fivebf
 \textfont1=\eightbfit \scriptfont1=\sixbfit \scriptscriptfont1=\fivebfit
 \textfont2=\eightbsy  \scriptfont2=\sixbsy  \scriptscriptfont2=\fivebsy
 \eightbf
 \baselineskip=10pt}

\def\tenpointbf{%
 \textfont0=\tenbf   \scriptfont0=\sevenbf   \scriptscriptfont0=\fivebf
 \textfont1=\tenbfit \scriptfont1=\sevenbfit \scriptscriptfont1=\fivebfit
 \textfont2=\tenbsy  \scriptfont2=\sevenbsy  \scriptscriptfont2=\fivebsy
 \tenbf}
        
\def\twelvepointbf{%
 \textfont0=\twelvebf   \scriptfont0=\eightbf   \scriptscriptfont0=\sixbf
 \textfont1=\twelvebfit \scriptfont1=\eightbfit \scriptscriptfont1=\sixbfit
 \textfont2=\twelvebsy  \scriptfont2=\eightbsy  \scriptscriptfont2=\sixbsy
 \twelvebf
 \baselineskip=14.4pt}

\def\seventeenpointbf{%
 \textfont0=\seventeenbf  \scriptfont0=\twelvebf  \scriptscriptfont0=\eightbf
 \textfont1=\seventeenbfit\scriptfont1=\twelvebfit\scriptscriptfont1=\eightbfit
 \textfont2=\seventeenbsy \scriptfont2=\twelvebsy \scriptscriptfont2=\eightbsy
 \seventeenbf
 \baselineskip=20.736pt}
 

\newdimen\srdim \srdim=\hsize
\newdimen\irdim \irdim=\hsize
\def\NOSECTREF#1{\noindent\hbox to \srdim{\null\dotfill ???(#1)}}
\def\SECTREF#1{\noindent\hbox to \srdim{\csname REF\romannumeral#1\endcsname}}
\def\INDREF#1{\noindent\hbox to \irdim{\csname IND\romannumeral#1\endcsname}}
\newlinechar=`\^^J
\def\openauxfile{
  \immediate\openin1\jobname.
  \ifeof1
  \let\sectref=\NOSECTREF \let\indref=\NOSECTREF
  \else
  \input \jobname.aux
  \let\sectref=\SECTREF \let\indref=\INDREF
  \fi
  \immediate\openout1=\jobname.aux
  \let\END=\end \def\end{\immediate\closeout1\END}}
        
\newbox\titlebox   \setbox\titlebox\hbox{\hfil}
\newbox\sectionbox \setbox\sectionbox\hbox{\hfil}
\def\folio{\ifnum\pageno=1 \hfil \else \ifodd\pageno
           \hfil {\eightpoint\copy\sectionbox\kern8mm\number\pageno}\else
           {\eightpoint\number\pageno\kern8mm\copy\titlebox}\hfil \fi\fi}
\footline={\hfil}
\headline={\folio}

\def\titlerunning#1{\setbox\titlebox\hbox{\eightpoint #1}}
\def\title#1{\noindent\hfil$\smash{\hbox{\seventeenpointbf #1}}$\hfil
             \titlerunning{#1}\medskip}

\newcount\numbersection \numbersection=-1
\def\sectionrunning#1{\setbox\sectionbox\hbox{\eightpoint #1}
  \immediate\write1{\string\def \string\REF 
      \romannumeral\numbersection \string{%
      \noexpand#1 \string\dotfill \space \number\pageno \string}}}
\def\section#1{%
  \par\vskip0.666cm\penalty -100
  \vbox{\baselineskip=14.4pt\noindent{{\twelvepointbf #1}}}
  \vskip2pt
  \penalty 500
  \advance\numbersection by 1
  \sectionrunning{#1}}

\def\subsection#1|{%
  \par\vskip0.5cm\penalty -100
  \vbox{\noindent{{\tenpointbf #1}}}
  \vskip1pt
  \penalty 500}

\newcount\numberindex \numberindex=0  
\def\index#1#2{%
  \advance\numberindex by 1
  \immediate\write1{\string\def \string\IND #1%
     \romannumeral\numberindex \string{%
     \noexpand#2 \string\dotfill \space \string\S \number\numbersection, 
     p.\string\ \space\number\pageno \string}}}

\newdimen\itemindent \itemindent=\parindent

\def\item#1{\par\noindent\hangindent\itemindent%
            \rlap{#1}\kern\itemindent\ignorespaces}
\def\itemitem#1{\par\noindent\hangindent2\itemindent%
            \kern\itemindent\rlap{#1}\kern\itemindent\ignorespaces}
\def\itemitemitem#1{\par\noindent\hangindent3\itemindent%
            \kern2\itemindent\rlap{#1}\kern\itemindent\ignorespaces}

\long\def\claim#1|#2\endclaim{\par\vskip 5pt\noindent 
{\tenpointbf #1.}\ {\sl #2}\par\vskip 5pt}

\def\proof{\noindent{\sl Proof}}

\def\today{\ifcase\month\or
January\or February\or March\or April\or May\or June\or July\or August\or
September\or October\or November\or December\fi \space\number\day,
\number\year}

\catcode`\@=11
\newcount\@tempcnta \newcount\@tempcntb 
\def\timeofday{{%
\@tempcnta=\time \divide\@tempcnta by 60 \@tempcntb=\@tempcnta
\multiply\@tempcntb by -60 \advance\@tempcntb by \time
\ifnum\@tempcntb > 9 \number\@tempcnta:\number\@tempcntb
  \else\number\@tempcnta:0\number\@tempcntb\fi}}
\catcode`\@=12

\def\bibitem#1&#2&#3&#4&%
{\hangindent=1.66cm\hangafter=1
\noindent\rlap{\hbox{\eightpointbf #1}}\kern1.66cm{\rm #2}{\sl #3}{\rm #4.}} 


\def\bC{{\Bbb C}}

\def\bN{{\Bbb N}}
\def\bP{{\Bbb P}}

\def\bR{{\Bbb R}}

\def\gm{{\euf m}}

\def\cE{{\Cal E}}
\def\cF{{\Cal F}}
\def\cO{{\Cal O}}
\def\cH{{\Cal H}}
\def\cI{{\Cal I}}
\def\cJ{{\Cal J}}
\def\cP{{\Cal P}}

\def\cZ{{\Cal Z}}


\def\square{{\hfill \hbox{
\vrule height 1.453ex  width 0.093ex  depth 0ex
\vrule height 1.5ex  width 1.3ex  depth -1.407ex\kern-0.1ex
\vrule height 1.453ex  width 0.093ex  depth 0ex\kern-1.35ex
\vrule height 0.093ex  width 1.3ex  depth 0ex}}}
\def\qed{\phantom{$\quad$}\hfill$\square$\medskip}
\def\hexnbr#1{\ifnum#1<10 \number#1\else
 \ifnum#1=10 A\else\ifnum#1=11 B\else\ifnum#1=12 C\else
 \ifnum#1=13 D\else\ifnum#1=14 E\else\ifnum#1=15 F\fi\fi\fi\fi\fi\fi\fi}
\def\msatype{\hexnbr\msafam}
\def\msbtype{\hexnbr\msbfam}
\mathchardef\restriction="3\msatype16   
\mathchardef\compact="3\msatype62
\mathchardef\smallsetminus="2\msbtype72   \let\ssm\smallsetminus
\mathchardef\subsetneq="3\msbtype28
\mathchardef\supsetneq="3\msbtype29
\mathchardef\leqslant="3\msatype36   \let\le\leqslant
\mathchardef\geqslant="3\msatype3E   \let\ge\geqslant
\mathchardef\ltimes="2\msbtype6E
\mathchardef\rtimes="2\msbtype6F

\let\ol=\overline

\let\text=\hbox
\def\build#1|#2|#3|{\mathrel{\mathop{\null#1}\limits^{#2}_{#3}}}

\def\Div{\mathop{\rm Div}\nolimits}
\def\lc{\mathop{\rm lc}\nolimits}

\def\Vol{\mathop{\rm Vol}\nolimits}

\def\dbar{{\overline\partial}}
\def\ddbar{{\partial\overline\partial}}


\def\loc{{\rm loc}}

\def\bfe{{\bf e}}

\long\def\InsertFig#1 #2 #3 #4\EndFig{\par
\hbox{\hskip #1mm$\vbox to#2mm{\vfil\special{" 
(/home/demailly/psinputs/grlib.ps) run
#3}}#4$}}
\long\def\LabelTeX#1 #2 #3\ELTX{\rlap{\kern#1mm\raise#2mm\hbox{#3}}}


\openauxfile

\title{Estimates on Monge-Amp\`ere}
\smallskip
\title{operators derived from}
\smallskip
\title{a local algebra inequality}
\titlerunning{Estimates on Monge-Amp\`ere operators}
\medskip
\centerline{\twelvebf Jean-Pierre Demailly}
\medskip
\centerline{Universit\'e de Grenoble I, D\'epartement de Math\'ematiques}
\centerline{Institut Fourier, 38402 Saint-Martin d'H\`eres, France}
\centerline{{\it e-mail\/}: {\tt demailly@fourier.ujf-grenoble.fr}}

\smallskip
\vskip7pt
\line{\hfill \it Dedicated to Professor Christer Kiselman, on the occasion
of his retirement} 
\vskip15pt

\noindent
{\bf Abstract.}
The goal of this short note is to relate the integrability property of
the exponential $e^{-2\varphi}$ of a plurisubharmonic function
$\varphi$ with isolated or compactly sup\-por\-ted singularities, to a
priori bounds for the Monge-Amp\`ere mass of $(dd^c\varphi)^n$. The
inequality is valid locally or globally on an arbitrary open subset 
$\Omega$ in~$\bC^n$. We show that $\int_\Omega(dd\varphi)^n<n^n$ implies
$\int_Ke^{-2\varphi}<+\infty$ for every compact subset~$K$ in~$\Omega$,
while functions of the form $\varphi(z)=n\log|z-z_0|$, $z_0\in\Omega$,
appear as limit cases. The result is derived from an inequality of pure 
local algebra, which turns out a posteriori to be equivalent to it, 
proved by A.~Corti in dimension $n=2$, and later extended by L.~Ein, 
T.~De Fernex and M.~Musta\c{t}\v{a} to arbitrary dimensions.
\medskip

\noindent
{\bf R\'esum\'e.}
Le but de cette note est d'\'etablir une relation entre la
propri\'et\'e d'int\'egrabilit\'e de l'exponentielle $e^{-2\varphi}$
d'une fonction plurisubharmonique $\varphi$ dont les singularit\'es
sont isol\'ees ou \`a support compact, et la donn\'ee de bornes a priori
pour la masse de Monge-Amp\`ere de $(dd^c\varphi)^n$.
L'in\'egalit\'e obtenue a lieu aussi bien localement que globalement,
ceci sur un ouvert arbitraire $\Omega$ de~$\bC^n$.  Nous montrons que
l'hypoth\`ese $\int_\Omega(dd\varphi)^n<n^n$ entra\^{\i}ne
$\int_Ke^{-2\varphi}<+\infty$ pour tout sous-ensemble compact~$K$
de~$\Omega$, les fonctions de la forme $\varphi(z)=n\log|z-z_0|$,
$z_0\in\Omega$, apparaissant comme des cas limites. Le r\'esultat se d\'eduit
d'une pure in\'egalit\'e d'alg\`ebre locale, qui se trouve a posteriori 
lui \^etre \'equivalente, successivement d\'emontr\'ee par A.~Corti 
en dimension $n=2$, puis \'etendue par L.~Ein, T.~De Fernex et
M.~Musta\c{t}\v{a} en dimensions arbitraires.
\bigskip

\noindent
{\bf Key words.}
Monge-Amp\`ere operator, local algebra, monomial ideal, 
Hilbert-Samuel multiplicity,
log-canonical threshold, plurisubharmonic function, Ohsawa-Takegoshi
$L^2$ extension theorem, approximation of singularities, birational rigidity.
\medskip
\noindent
{\bf Mots-cl\'es.}
Op\'erateur de Monge-Amp\`ere, alg\`ebre locale, id\'eal monomial,
multiplicit\'e de Hilbert-Samuel, seuil log-canonique, fonction
plurisousharmonique, th\'eor\`eme d'extension $L^2$ d'Ohsawa-Takegoshi,
approximation des singularit\'es, rigidit\'e birationnelle.
\bigskip

\noindent
{\bf AMS Classification.} 32F07, 14B05, 14C17

\section{1. Main result}

Here we put $d^c={i\over 2\pi}(\dbar-\partial)$ so that
$dd^c={i\over\pi}\ddbar$. The normalization of the $d^c$ operator is
chosen such that we have precisely $(dd^c\log|z|)^n=\delta_0$ for the
Monge-Amp\`ere operator in $\bC^n$. The Monge-Amp\`ere operator is
defined on locally bounded plurisubharmonic functions according to the
definition of Bedford-Taylor [BT76, BT82]; it can also be extended to
plurisubharmonic functions with isolated or compactly supported poles
by [Dem93]. Our main result is the following a priori estimate for the
Monge-Amp\`ere operator acting on functions with compactly supported
poles.

\claim (1.1) Main Theorem|Let $\Omega$ be an open subset 
in $\bC^n$, $K$~a compact subset of~$\Omega$, and let $\varphi$ be a 
plurisubharmonic function on $\Omega$ such that $-A\le\varphi\le 0$
on $\Omega\ssm K$ and
$$
\int_\Omega(dd^c\varphi)^n\le M<n^n.
$$
Then there is an a priori upper bound for the Lebesgue integral of
$e^{-2\varphi}$, namely
$$
\int_K e^{-2\varphi}d\lambda\le C(\Omega,K,A,M),
$$
where the constant $C(\Omega,K,A,M)$ depends on the given parameters
but not on the function~$\varphi$.
\endclaim

We first make a number of elementary remarks.
\medskip

\noindent {\bf (1.2)} The result is optimal as far as the
Monge-Amp\`ere bound $M<n^n$ is concerned, since functions
$\varphi_\varepsilon(z)=(n-\varepsilon)\log|z-z_0|$, $z_0\in K^\circ
\subset \Omega$ satisfy
$\int_\Omega(dd^c\varphi_\varepsilon)^n=(n-\varepsilon)^n$, but
$\int_Ke^{-2\varphi_\varepsilon}d\lambda$ tends to $+\infty$ as
$\varepsilon$ tends to zero.\medskip

\noindent {\bf (1.3)} The assumption $-A\le\varphi\le 0$ on
$\Omega\ssm K$ is required, as it forces the poles of $\varphi$ to be
compactly supported -- a condition needed to define properly the
Monge-Amp\`ere measure $(dd^c\varphi)^n$ (see e.g.\ [Dem93]). In 
any case, the functions
$\varphi_\varepsilon(z)={1\over 2}\ln(|z_1|^2+\varepsilon^2)$ satisfy
$\int_\Omega(dd^c\varphi_\varepsilon)^n=0<n^n$, but $\int_K
e^{-2\varphi_\varepsilon}d\lambda$ is unbounded as $\varepsilon$ tends
to $0$,  whenever $K$ contains at least one interior point located on the
hyperplane $z_1=0$. The limit $\varphi(z)=\ln|z_1|$ of course does not
have compactly supported poles. In such a circumstance, C.O.~Kiselman [Ki84]
observed long ago that the Monge-Amp\`ere mass of $(dd^c\varphi)^n$ need
not be finite or well defined.
\medskip

\noindent {\bf (1.4)} The a priori estimate (1.1) can be seen as a non linear
analogue of Skoda's criterion for the local integrability of $e^{-2\varphi}$.
Let us recall Skoda's criterion$\,$:
{\it if~the Lelong number $\nu(\varphi,z_0)$ satisfies $\nu(\varphi,z_0)<1$, 
then $e^{-2\varphi}$ is locally integrable near~$z_0$, and if 
 $\nu(\varphi,z_0)\ge n$, then $\int_Ve^{-2\varphi}d\lambda=+\infty$ on every
neighborhood $V$ of~$z_0$}. The gap between $1$ and $n$ is an important
feature of potential theory in several complex variables, and it therefore
looks like an interesting bonus that there is no similar discrepancy for 
the estimate given by Theorem~1.1. One of the reasons is that 
$(dd^c\varphi)^n$ takes into account all dimensions simultaneously, while
the Lelong number only describes the minimal vanishing order with respect to
arbitrary lines (or holomorphic curves).\qed
\medskip

The proof consists of several steps, the main of which is a reduction to the 
following result of local algebra, due to A.~Corti [Cor00] in dimension
$2$ and L.~Ein, T.~De Fernex and M.~Musta\c{t}\v{a} [dFEM04] in general.

\claim (1.5) Theorem|Let $\cJ$ be an ideal in the rings of germs 
$\cO_{\bC^n,0}$ of holomorphic functions in $n$ variables, such that
the zero variety $V(\cJ)$ consists of the single point~$\{0\}$.
Let $e(\cJ)$ be the Samuel multiplicity of $\cJ$, i.e.\
$$
e(\cJ)=\lim_{k\to+\infty}{n!\over k^n}\dim\cO_{\bC^n,0}/\cJ^k,
$$
the leading coefficient in the Hilbert polynomial of~$\cJ$.
Then the log canonical threshold of $\cJ$ satisfies
$$
\lc(\cJ)\ge {n \over e(\cJ)^{1/n}},
$$
and the equality case occurs if and only if
the integral closure $\ol\cJ$ of $\cJ$ is a power $\gm^s$ of 
the maximal ideal.
\endclaim 

Recall that the {\it log canonical threshold} $\lc(\cJ)$ of an ideal $\cJ$
is the supremum of all numbers $c>0$ such that
$(|g_1|^2+\ldots+|g_N|^2)^{-c}$ is integrable near $0$ for any set of
generators $g_1,\ldots,g_N$ of $\cJ$. If this supremum is less than
$1$ $\big[$this is always the case after replacing $\cJ$ by a sufficiently
high power $\cJ^m$, which yields $e(\cJ^m)=m^ne(\cJ)$
and $\lc(\cJ^m)={1\over m}\lc(\cJ)\,\big]$, the integrability condition
exactly means that
the divisor $cD$ associated with a generic element $D=\Div_f$, $f\in\cJ$, is
{\it Kawamata log terminal} (klt), i.e.\ that after blowing up and resolving
the singularities to get a divisor with normal crossings, the
associated divisor $\mu^*(cD)-E$ has coefficients${}<1$, where $\mu$ is the 
blow-up map and $E$ its jacobian divisor (see e.g.\ [DK00] for details).

In fact, Theorem 1.5 follows from Theorem 1.1 by taking $\Omega$
equal to a small ball $B(0,r)\subset\bC^n$ and
$\varphi(z)={c\over 2}\log\sum_j|g_j|^2$ where $g_1,\ldots,g_N$
are local generators of $\cJ$. For this, one observes that the Monge-Amp\`ere
mass of $(dd^c\varphi)^n$ carried by $\{0\}$ is equal to
$c^ne(\cJ)$ (Lemma 2.1 below), hence the integrability of 
$(|g_1|^2+\ldots+|g_N|^2)^{-c}=e^{-2\varphi}$ holds true as soon
as $c^ne(\cJ)<n^n\,$; notice that the integral $\int_{B(0,r)}(dd^c\varphi)^n$
converges to the mass carried by $0$ as the radius $r$ tends to zero.

However, the strategy of the proof goes the other way round$\,$: Theorem 1.1
will actually be derived from Theorem 1.5 by means of the approximation
techniques for plurisubharmonic functions developped in [Dem92] and
the result on semi-continuity of singularity exponents ${}={}$log canonical
thresholds) obtained in [DK00]. It is somewhat strange that one has to make
a big detour through local algebra (and approximation of analytic objects
by polynomials, as in [DK00]), to prove what finally appears to be
a pure analytic estimate on Monge-Amp\`ere operators. 

It would be interesting to know whether a direct proof can be obtained
by methods which are more familiar to analysts (integration by parts,
convexity inequalities, integral kernels for $\dbar\;\ldots$). One of
the consequences of our use of a ``purely qualitative'' algebraic detour 
is that the constants $C(\Omega,K,A,M)$ appearing in Theorem~1.1 are non
effective. On the other hand, we would like to know what kind of
dependance this constants have e.g.\ on $n^n-M>0$, and also what are the
extremal functions (for instance in the case when $\Omega$ and $K$ are
concentric balls). The question is perhaps more difficult than it would
first appear, since the most obvious guess is that the extremal functions
are singular ones with a logarithmic pole
$\varphi(z)\sim\lambda\log|z-z_0|\,$; the reason for this expectation
is that the equality case in Theorem 1.5 is achieved precisely when
the integral closure of the ideal $\cJ$ is equal to a power of the
maximal ideal.

I would like to thank I.~Chel'tsov, R.~Lazarsfeld, L.~Ein, T.~de
Fernex, M.~Musta\c{t}\v{a} for explaining to me the algebraic issues
involved in the inequalities just discussed (see e.g.\ [Che05]). It is
worth mentioning that inequality 1.5 is related to deep questions of
algebraic geometry such as the birational (super)rigidity of Fano
manifolds$\,$; for instance, following ideas of Corti and Pukhlikov
([Cor95], [Cor00]), it is proved in [dFEM03] that every smooth
hypersurface of degree $N$ in $\bP^N$ is birationally superrigid at
least for $4\le N\le 12$, hence that such a hypersurface cannot be
rational -- this is a far reaching generalization of the classical
result by Iskovskikh-Manin ([IM72], [Isk01]) that $3$-dimensional
quartics are not rational.

I am glad to dedicate this paper to Professor C.O.~Kiselman whose work
has been a great source of inspiration for my own research in complex
analysis, especially on all subjects related to Monge-Amp\`ere
operators, Lelong numbers and attenuation of singularities of
plurisubharmonic functions ([Kis78, 79, 84, 94a, 94b]). Various
incarnations of these concepts and results appear throughout the
present paper.
\bigskip

\section{2. Proof of the integral inequality}

The first step is to related Monge-Amp\`ere masses to Samuel multiplicities.
The relevant result is probably known, but we have not been able to find 
a precise reference in the litterature.

\claim (2.1) Lemma|In a neighborhood of $0\in\bC^n$, let
$\varphi(z)={1\over 2}\log\sum_{j=1}^N|g_j|^2$ where $g_1,\ldots,g_N$
are germs of holomorphic functions which have $0$
as their only common zero. Then the Monge-Amp\`ere
mass of $(dd^c\varphi)^n$ carried by $\{0\}$ is equal to the
Samuel multiplicity $e(\cJ)$ of the
ideal $\cJ=(g_1,\ldots,g_N)\subset\cO_{\bC^n,0}$.
\endclaim

\proof. For any point $a\in G_{N,n}$, the Grassmannian of $n$-dimensional
subspaces in~$\bC^N$, we define
$$
\varphi_a(z)={1\over 2}
\log\sum_{i=1}^n\Big|\sum_{j=1}^N\lambda_{ij}g_j(z)\Big|^2
$$
where $(\lambda_{jk})$ is the $n\times N$-matrix of an orthonormal basis
of the subspace~$a$. It is easily shown that $\varphi_a$ is defined in a 
unique way and that we have the Crofton type formula
$$
(dd^c\varphi(z))^n=\int_{a\in G_{N,n}} (dd^c\varphi_a(z))^nd\mu(a)\leqno(2.2)
$$
where $\mu$ is the unique $U(N)$-invariant probability measure on the
Grassmannian. In fact this can be proved from the related equality
$$
(dd^c\log|w|^2)^n=\int_{a\in G_{N,n}} (dd^c\log|\pi_a(w)|^2)^nd\mu(a)
$$
where $\bC^N\ni w\mapsto\pi_a(w)$ is the orthogonal projection onto 
$a\subset\bC^n$, which itself follows by unitary invariance and a
degree argument (both sides have degree one as bidegree $(n,n)$ currents
on the projective space $\bP^{N-1}$). One then applies the 
substitution $w=g(z)$ to get the 
general case. The right hand side of (2.2) is well defined since the
poles of $\varphi_a$ form a finite set for a generic point $a$ in the 
grassmannian$\,$;
then $(dd^c\varphi_a(z))^n$ is just a sum of Dirac masses with integral
coefficients (the local degree of the corresponding germ of map $g_a:
z\mapsto \big(\sum_{1\le j\le N}\lambda_{ij}g_j(z)\big)_{1\le i\le N}$ from 
$\bC^n$ to $\bC^n$ near the given point). By a continuity argument, the 
coefficient of $\delta_0$ is constant except on some analytic stratum
in the Grasmannian, and by Fubini, the mass carried by $(dd^c\varphi)^n$
at $0$ is thus equal to the degree of $g_a$ at $0$ for generic~$a$. Now, it 
is a well-known fact of commutative algebra that the Hilbert-Samuel
multiplicity $e(\cJ)$ is equal to the intersection number of the divisors
associated with a generic $n$-tuple of elements of~$\cJ$ (Bourbaki,
Alg\`ebre Commutative [BAC83], VIII 7.5, Prop.~7). That intersection number
is also equal to the generic value of the Monge-Amp\`ere mass 
$$
\big(dd^c\log|\lambda_1\cdot g|\big)\wedge\ldots\wedge
\big(dd^c\log|\lambda_n\cdot g|\big)(0).
$$
By averaging with respect to the $\lambda_j\,$'s, this appears to be the 
same as the generic value of~$(dd^c\varphi_a)^n(0)$.\qed
\medskip

We now briefly recall the ideas involved in the proof of Theorem~1.5, as taken 
from [dFEM04]. In order to prove the main inequality of 1.5 (which can be
rewritten as $e(\cJ)\ge n^n/\lc(\cJ)^n$), it is sufficient
so show that
$$
\dim\cO_{\bC^n, 0}/\cJ \ge n^n/(n!\lc(\cJ)^n).\leqno(2.3)
$$
In fact, since by definition $\lc(\cJ^k)={1\over k}\lc(\cJ)$, a substitution
of $\cJ$ by $\cJ^k$ in (2.3) yields
$$
{n!\over k^n}\dim\cO_{\bC^n, 0}/\cJ^k \ge n^n/\lc(\cJ)^n
$$
an we get the expected inequality 1.5 by letting $k$ tend to $+\infty$. 
Now, by fixing a multiplicative order on the coordinates $z_j$, it is 
well known that one can construct a flat family $(\cJ_s)_{s\in\bC,0}$
depending on a small complex parameter~$s$ such that $\cJ_0$ is a monomial
ideal and $\cO_{\bC^n,0}/\cJ_s\simeq\cO_{\bC^n,0}/\cJ$ for all $s\ne 0$
(see Eisenbud [Ei95] for a nice discussion in the algebraic case).
The semicontinuity property of the log canonical threshold (see for example
[DK00]) implies that $\lc(\cJ_0)\le\lc(\cJ_s)$ for small~$s$. 

The proof is then reduced to the case when $\cJ$ is a monomial ideal, i.e.\
an ideal generated by a family of monomials $(z^{\beta_j})$.
In the latter situation, the argument proceeds from an explicit formula
for $\lc(\cJ)$ due to J.~Howald [Ho01]~: let $P(\cJ)$ be the Newton polytope
of $\cJ$, i.e.\ the convex hull of the points $\beta\in\bN^n$ associated with
all monomials $z^\beta\in\cJ\,$; then putting $\bfe=(1,\ldots,1)\in\bN^n$,
$$
{1\over\lc(\cJ)}=\min\big\{\alpha>0\,;\;\alpha\cdot\bfe\in P(\cJ)\big\}
$$
(the reader can take this as a clever exercise on the convergence of
integrals defined by sums of monomials in the denominator).  Let $F$ be the
facet of $P(\cJ)$ which contains the point ${1\over\lc(\cJ)}\bfe$, and let
$\sum x_j/a_j=1$, $a_j>0$ be the equation of this hypersurface in $\bR^n$.
Let us denote also by $F_+$ and $F_-$ the open half-spaces delimited by $F$,
such that $\bR_+^n\cap F_-$ is relatively compact and $\bR_+^n\cap F_+$ is 
unbounded. Then $\Vol(\bR_+^n\cap F_-)={1\over n!}\prod a_j$ and therefore, 
since 
$\bR_+^n\ssm P(\cJ)$ contains $\bR_+^n\cap F_-$, we get
$$
\Vol(\bR_+^n\ssm P(\cJ))\ge \Vol(\bR_+^n\cap F_-)= {1\over n!}\prod a_j.
$$
On the other hand, $\dim\cO_{\bC^n,0}/\cJ$ is at least equal to the number 
of elements of $\bN^n\ssm P(\cJ)$, which is itself at least equal to
$\Vol(\bR_+^n\ssm P(\cJ))$ since the unit cubes $\beta+[0,1]^n$ with
$\beta\in \bN^n\ssm P(\cJ)$ cover the complement $\bR_+^n\ssm P(\cJ)$.
This yields
$$
\dim\cO_{\bC^n,0}/\cJ\ge {1\over n!}\prod a_j.
$$
As ${1\over\lc(\cJ)}\bfe$ belongs to $F$, we have $\sum 1/a_j=\lc(\cJ)$.
The inequality between geometric and arithmetic means implies
$$
\Big(\prod {1\over a_j}\Big)^{1/n}\le {1\over n}\sum{1\over a_j} =
{\lc(\cJ)\over n}
$$
and inequality (2.3) follows.
We refer to [dFEM04] for the discussion of the equality case.\qed
\medskip

The next ingredient is the following basic approximation theorem for 
plurisubharmonic functions through the Bergman kernel trick and the 
Ohsawa-Takegoshi theorem [OT87], the first version of which appeared
in~[Dem92]. We start with the general concept of complex singularity
exponent introduced in [DK00], which extends the concept of log 
canonical threshold.

\claim (2.4) Definition|Let $X$ be a complex manifold and $\varphi$ be
a plurisubharmonic $($psh$)$ function on~$X$. For any compact set
$K\subset X$, we introduce the ``complex singularity exponent'' of
$\varphi$ on $K$ to be the nonnegative number
\index{CO}{Complex singularity exponent}
$$
c_K(\varphi)=\sup\big\{c\ge 0\,:\,\exp(-2c\varphi)
\text{ is $L^1$ on a neighborhood of $K$}\big\},
$$
and we define the ``Arnold multiplicity'' to be
$\lambda_K(\varphi)=c_K(\varphi)^{-1}\,:$
$$
\lambda_K(\varphi)=\inf\big\{\lambda>0\,:\,\exp(-2\lambda^{-1}\varphi)
\text{ is $L^1$ on a neighborhood of $K$}\big\}.
$$
\endclaim

In the case where $\varphi(z)={1\over 2}\log\sum |g_j|^2$, the exponent
$c_{z_0}(\varphi)$ is the same as the log canonical threshold of
the ideal $\cJ=(g_j)$ at the point~$z_0$.

\claim (2.5) Theorem {\rm([Dem92, DK00])}|Let $\varphi$ be a plurisubharmonic 
function on a
bounded pseudoconvex open set~$\Omega\subset\bC^n$. For every real 
number $m>0$, let
$\cH_{m\varphi}(\Omega)$ be the Hilbert space of holomorphic functions $f$
on $\Omega$ such that $\int_\Omega|f|^2e^{-2m\varphi}dV<+\infty$ and
let $\psi_m={1\over 2m}\log\sum_k|g_{m,k}|^2$ where $(g_{m,k})_k$
is an orthonormal basis of~$\cH_{m\varphi}(\Omega)$. Then\/$:$
\smallskip
\item{\rm (a)} There are constants $C_1,C_2>0$ independent of $m$ and
$\varphi$ such that
$$
\varphi(z)-{C_1\over m}\le
\psi_m(z)\le\sup_{|\zeta-z|<r}\varphi(\zeta)+{1\over m}\log{C_2\over r^n}
$$
for every $z\in\Omega$ and $r<d(z,\partial\Omega)$. In particular,
$\psi_m$ converges to $\varphi$ pointwise and in $L^1_{\rm loc}$ topology
on~$\Omega$ when $m\to+\infty$ and
\smallskip
\item{\rm (b)} The Lelong numbers of $\varphi$ and $\psi_m$ are related by
$$
\nu(\varphi,z)-{n\over m}\le\nu(\psi_m,z)\le
\nu(\varphi,z)\quad\text{ for every $z\in\Omega$.}
$$
\item{\rm (c)} For every compact set $K\subset\Omega$, the Arnold 
multiplicity of $\varphi$,~$\psi_m$ and of the multiplier ideal sheaves
$\cI(m\varphi)$ are related by
$$
\lambda_K(\varphi)-{1\over m}\le\lambda_K(\psi_m)=
{1\over m}\lambda_K(\cI(m\varphi))\le \lambda_K(\varphi).
$$
\vskip0pt
\endclaim 

The final ingredient is the following fundamental semicontinuity result 
from [DK00].

\claim (2.6) Theorem {\rm([DK00])}|Let $X$ be a complex manifold. 
Let $\cZ^{1,1}_+(X)$
denote the space of closed positive currents of type $(1,1)$ on~$X$, equipped
with the weak topology, and let $\cP(X)$ be the set of locally $L^1$ psh 
functions on $X$, equipped with the topology of $L^1$ convergence on 
compact subsets $(=$~topology induced by the weak topology$)$. Then
\smallskip
\item{\rm(a)} The map $\varphi\mapsto c_K(\varphi)$ is lower semi-continuous 
on $\cP(X)$, and the map $T\mapsto c_K(T)$ is lower semi-continuous on
$\cZ^{1,1}_+(X)$.
\smallskip
\item{\rm(b)} {\rm(``Effective version'')}. Let $\varphi\in\cP(X)$ be given.
If $c<c_K(\varphi)$ and $\psi$ converges to $\varphi$ in $\cP(X)$, then
$e^{-2c\psi}$ converges to $e^{-2c\varphi}$ in $L^1$ norm over some
neighborhood $U$ of~$K$.
\vskip0pt
\endclaim

\noindent {\bf (2.7) Proof of Theorem 1.1.} Assume that the conclusion
of theorem 1.1 is wrong. Then there exist a compact set $K\subset \Omega$,
constants $M<n^n$, $A>0$ and a sequence $\varphi_j$ of plurisubharmonic
functions such that $-A\le\varphi_j\le 0$ on $\Omega\ssm K$ and
$\int_\Omega(dd^c\varphi_j)^n\le M$, while $\int_K e^{-2\varphi_j}d\lambda$ 
tends to $+\infty$ as $j$ tends to $+\infty$. By~well-known properties of
potential theory, the condition $-A\le\varphi_j\le 0$ on
$\Omega\ssm K$ ensures that the
sequence $(\varphi_j)$ is relatively compact in the $L^1_\loc$ topology
on $\Omega\,$:
in fact, the Laplacian $\Delta\varphi_j$ is a uniformly bounded measure
on every compact of $\Omega\ssm K$, and this property extends to all compact
subsets of $\Omega$ by Stokes'
theorem and the fact that there is a strictly subharmonic function 
on $\Omega\,$; we then conclude by an elementary (local) Green kernel 
argument. Therefore there exists a subsequence of $(\varphi_j)$
which converges almost everywhere and
in $L^1_\loc$ topology to a limit $\varphi$ such that $-A\le\varphi\le 0$
on $\Omega\ssm K$ and $\int_\Omega(dd^c\varphi)^n\le M$. On the other 
hand, we must have $c_K(\varphi)\le 1$ by (2.6$\,$b)
(hence $\int_Ke^{-2(1+\varepsilon)\varphi}d\lambda=+\infty$ for every
$\varepsilon>0$). 

As $c_K(\varphi)=\inf_{z\in K}c_{\{z\}}(\varphi)$ and
$z\mapsto c_{\{z\}}(\varphi)$ is lower semicontinuous, there exists
a point $z_0\in K$ such that $c_{\{z_0\}}(\varphi)\le 1$.
By theorem 2.5 applied on a small ball $B(z_0,r)$, we can 
approximate $\varphi$ by a sequence
of psh functions of the form $\psi_m={1\over 2m}\log\sum|g_{m,k}|^2$
on $B(z_0,r)$. Inequality (2.5$\,$c) shows that we have 
$$
c_{\{z_0\}}(\psi_m)\le{1\over 1/c_{\{z_0\}(\varphi)}-1/m}\le{1\over 1-1/m}
<1+\varepsilon
$$
for $m$
large, hence $c_{\{z_0\}}((1+\varepsilon)\psi_m)\le 1$.
However, the analytic strata of positive Lelong numbers
of $\varphi$ must be contained in~$K$, hence they are isolated points 
in~$\Omega$, and thus the poles of $\psi_m$ are isolated. By the weak 
continuity of the Monge-Amp\`ere operator, we have
$$\int_{\ol B(z_0,r')}(dd^c(1+\varepsilon)\psi_m)^n\le(1+\varepsilon)^{n+1}
\int_{\ol B(z_0,r')}(dd^c\varphi)^n\le (1+\varepsilon)^{n+1}M^n$$
for $m$ large, for any $r'<r$. If $\varepsilon$ is chosen so small that 
$(1+\varepsilon)^{n+1}M^n<n^n$, then
the Monge-Amp\`ere mass of $(1+\varepsilon)\psi_m$ at $z_0$
is strictly less than $n^n$, but the log canonical threshold is at most 
equal to~$1$. This contradicts inequality 1.5, when using Lemma 2.1 to
identify the Monge-Amp\`ere mass with the Samuel multiplicity.\qed

\claim (2.8) Remark|{\rm As the proof shows, the arguments are
mostly of a local nature (the main problem is to ensure convergence of
the integral of $e^{-2\psi}$ on a neighborhood of the poles of
an approximation $\psi$ of $\varphi$ with logarithmic poles). Therefore 
Theorem 1.1 is also valid for a plurisubharmonic 
function $\varphi$ on an arbitrary non singular complex variety $X$, 
provided that $X$ does not possess positive dimensional complex analytic 
subsets (any open subset $\Omega$ in a Stein manifold will thus do).
We leave the reader complete the obvious details.}
\endclaim

\claim (2.9) Remark|{\rm The proof is highly non constructive, so it seems
at this point that there is no way of producing an explicit bound
$C(\Omega,K,A,M)$. It would be interesting to find a method to calculate
such a bound, even a suboptimal one.}
\endclaim

\claim (2.10) Remark|{\rm The equality case in Theorem 1.5 suggests that
extremal functions with respect to the integral
$\int_Ke^{-2\varphi}d\lambda$ might be fonctions with Monge-Amp\`ere
measure $(dd^c\varphi)^n$ concentrated at one point $z_0\in K$, and a
logarithmic pole at~$z_0$. We are unsure what the correct boundary
conditions should be, so as to actually get nice extremal functions of
this form. We expect that an adequate condition is to assume that $\varphi$
has zero boundary values. Further potential theoretic arguments would
be needed for this, since prescribing the boundary values is not enough
to get the relative compactness of the family in the weak topology (but this
might be the case with the granted additional upper bound $n^n$ on the 
Monge-Amp\`ere mass){\parindent=0mm\footnote{${}^{(1)}$}{\eightpoint
\leftskip=5.3mm
After the present paper was completed,
Ahmed Zeriahi sent us a short proof of this fact, and also derived
a stronger integral bound valid on the whole of~$\Omega$. See the 
Appendix below.\vskip0pt}.}}
\endclaim

We end this discussion by stating two generalizations of theorem 1.1
whose algebraic counterparts are useful as well for their 
applications to algebraic geometry  (see [Che05] and [dFEM03]).

\claim (2.11) Theorem|Let $\Omega$ be an open subset in $\bC^n$, 
$K$~a compact subset of~$\Omega$, and let $\varphi, \psi$ be plurisubharmonic
functions on $\Omega$ such that $-A\le\varphi,\psi\le 0$ on $\Omega\ssm K$,
with $c_K(\psi)\ge {1\over \gamma}$, $\gamma<1$ and
$$
\int_\Omega(dd^c\varphi)^n\le M<n^n(1-\gamma)^n.
$$
Then 
$$
\int_K e^{-2\varphi-2\psi}d\lambda\le C(\Omega,K,A,\gamma,M),
$$
where the constant $C(\Omega,K,A,\gamma,M)$ depends on the 
given parameters but not on the functions~$\varphi$, $\psi$.
\endclaim

\noindent
{\it Proof.} This is an immediate consequence of H\"older's inequality
for the conjugate exponents $p=1/(1-\gamma-\varepsilon)$ and 
$q=1/(\gamma+\varepsilon)$, applied
to the functions $f=\exp(-2\varphi)$ and $g=\exp(-2\psi)\,$: 
when $\varepsilon>0$ is small enough,
the Monge-Amp\`ere hypothesis for $\varphi$ precisely implies that $f$ 
is in $L^p(K)$ thanks to Theorem~1.1, and the assumption $c_K(\psi)\ge 
{1\over\gamma}$ implies by definition that $g$ is in $L^q(K)$.\qed
\medskip

In the case where $D=\sum \gamma_j D_j$ is an effective divisor with normal
crossings and $\psi$ has codimension $1$ analytic singularities given
by $D$ (i.e.\ $\psi(z)\sim \sum\gamma_j\log|z_j|$ in suitable
local analytic coordinates), we see that Theorem 2.11 can be applied
with $\gamma=\max(\gamma_j)$ and with the Monge-Amp\`ere upper bound
$n^n(1-\max(\gamma_j))^n$. In~this circumstance, it turns out that the 
latter bound can be improved.

\claim (2.12) Theorem|Let $\Omega$ be an open subset in $\bC^n$, 
$K$~a compact subset of~$\Omega$, and $\varphi$ a plurisubharmonic
function on $\Omega$ such that $-A\le\varphi\le 0$ on $\Omega\ssm K$.
Assume that there are constants $0\le \gamma_1,\ldots,\gamma_n<1$ such that
$$
\int_\Omega(dd^c\varphi)^n\le M<n^n\prod_{1\le j\le n}(1-\gamma_j).
$$
Then 
$$
\int_K e^{-2\varphi(z)}\prod_{1\le j\le n}|z_j|^{-2\gamma_j}d\lambda
\le C(\Omega,K,A,\gamma_j,M),
$$
where the constant $C(\Omega,K,A,\gamma_j,M)$ depends on the 
given parameters but not on the function~$\varphi$.
\endclaim

\noindent
{\it Proof.} In the case when $\gamma_j$ is the form
$\gamma_j=1-1/p_j$ and $p_j\ge 1$ is an integer, Theorem 2.12 can
be derived directly from the arguments of the proof of Theorem~1.1. 
Since the estimate is
essentially local, we only have to check convergence near the poles
of $\varphi$, in the case when $\varphi$ has an isolated analytic pole 
located on the support of the divisor $D$. Assume that the pole is the 
center of a polydisk $D(0,r)=\prod
D(0,r_j)$, in coordinates chosen so that the components of $D$ are the
coordinates hyperplanes $z_j=0$. We simply apply Theorem 1.1 to
the function $\widetilde\varphi(z)=\varphi(z_1^{p_1},\ldots,
z_n^{p_n})$ (with $p_j=1$ if the component $z_j=0$ does not occur in
$D$).  We then get
$$
\int_{\prod D(0,r_j^{1/p_j})}(dd^c\widetilde\varphi)^n=
p_1\ldots p_n\int_{D(0,r)}(dd^c\varphi)^n=
\prod(1-\gamma_j)^{-1}\int_{D(0,r)}(dd^c\varphi)^n
$$
by a covering degree argument, while
$$
\int_{\prod\ol D(0,\rho_j^{1/p_j})}e^{-2\widetilde\varphi}d\lambda=
\int_{D(0,\rho)}e^{-2\varphi}\Big(\prod|z_j|^{2(1-1/p_j)}\Big)^{-1}
d\lambda
$$
by a change of variable $\zeta=z_j^{p_j}$. We do not have such a simple 
argument when the $\gamma_j\,$'s are arbitrary real numbers less than~$1$. 
In~that case, the proof consists of repeating the steps of Theorem~1.1,
with the additional observation that the statement of local 
algebra corresponding to Theorem 2.12 (i.e.\ with
$\varphi(z)=c\log\sum|g_j|^2$ possessing one isolated pole) is 
still valid by [dFEM03], Lemma~2.4.\qed
\bigskip\null

\section{References}
\medskip

{\eightpoint

\bibitem[BAC83]&Bourbaki, N.:& Alg\`ebre Commutative, chapter VIII et IX;&
Masson, Paris, 1983&
  
\bibitem[BT76]&Bedford, E., Taylor, B.A.:& The Dirichlet problem for a complex
Monge-Amp\`ere equation;& Invent.\ Math.\ {\bf 37} (1976) 1--44&

\bibitem[BT82]&Bedford, E., Taylor, B.A.:& A new capacity for plurisubharmonic
functions;& Acta Math.\ {\bf 149} (1982) 1--41&

\bibitem[Che05]&Chel'tsov, I.:& Birationally rigid Fano manifolds;&
  Uspekhi Mat.\ Nauk {\bf 60:}5 (2005), 71--160 and Russian Math.\ Surveys
  {\bf 60:}5 (2005), 875--965&

\bibitem[Cor95]&Corti, A.:&Factoring birational maps of threefolds after 
  Sarkisov;& J.\ Algebraic Geom.\ {\bf 4} (1995), 223--254&

\bibitem[Cor00]&Corti, A.:& Singularities of linear systems and $3$-fold
  birational geometry;& London Math.\ Soc.\ Lecture Note Ser.\ { \bf 281}
  (2000) 259--312&

\bibitem[dFEM03]&de Fernex, T., Ein, L., and Musta\c{t}\v{a}:& 
  Bounds for log canonical thresholds with applications to birational 
  rigidity;& Math.\ Res.\ Lett.\ {\bf 10} (2003) 219--236&

\bibitem[dFEM04]&de Fernex, T., Ein, L., and Musta\c{t}\v{a}:& 
  Multiplicities and log canonical thresholds;& J.\ Algebraic Geom.\
  {\bf 13} (2004) 603--615&

\bibitem[Dem90]&Demailly, J.-P.:& Singular hermitian metrics on
  positive line bundles;& Proceedings of the Bayreuth conference
  ``Complex algebraic varieties'', April~2-6, 1990, edited by
  K.~Hulek, T.~Peternell, M.~Schneider, F.~Schreyer, Lecture Notes in
  Math.\ n${}^\circ\,$1507, Springer-Verlag, 1992&

\bibitem[Dem92]&Demailly, J.-P.:& Regularization of closed positive currents 
  and Intersection Theory;& J.\ Alg.\ Geom.\ {\bf 1} (1992), 361--409&

\bibitem[Dem93]&Demailly, J.-P.:& Monge-Amp\`ere operators, Lelong numbers 
  and intersection theory;& Complex Analysis and Geometry, Univ.\ Series in 
  Math., edited by V.~Ancona and A.~Silva, Plenum Press, New-York, 1993&

\bibitem[DK00]&Demailly, J.-P., Koll\'ar, J.:& Semi-continuity of complex
singularity exponents and K\"ahler-Einstein metrics on Fano orbifolds;&
Ann.\ Sci.\ Ecole Norm.\ Sup.\ (4) {\bf 34} (2001), 525--556&

\bibitem[Ho01]&Howald, J.:& Multiplier ideals of monomial ideals;& Trans.\
  Amer.\ Math.\ Soc.\ {\bf 353} (2001), 2665--2671&

\bibitem[IM72]&Iskovskikh, V.~A.\  and Manin, Yu. I.:& Three-dimensional
  quartics and counterexamples to the L\"uroth problem;& Mat.\ Sb.\ 
  {\bf 86} (1971), 140--166; English transl., Math.\ Sb.\ {\bf 15} (1972), 
  141--166&

\bibitem[Isk01]&Iskovskikh, V.~A.:& Birational rigidity and Mori theory;&
  Uspekhi Mat.\ Nauk {\bf 56}:2 (2001) 3--86;
  English transl., Russian Math.\ Surveys {\bf 56}:2 (2001), 207--291.&

\bibitem[Kis78]&Kiselman, C.~O.:& The partial Legendre transformation
  for plurisubharmonic functions;& Inventiones Math.\ {\bf 49} (1978)
  137--148&

\bibitem[Kis79]&Kiselman, C.~O.:& Densit\'e des fonctions
  plurisousharmoniques;& Bulletin de la Soci\'et\'e Math\'ematique de
  France, {\bf 107} (1979) 295--304&

\bibitem[Kis84]&Kiselman, C.~O.:& Sur la d\'efinition de l'op\'erateur
  de Monge-Amp\`ere complexe;& Analyse Complexe; Proceedings of the
  Journ\'ees Fermat -- Journ\'ees SMF, Toulouse 1983$\,$; Lecture Notes in
  Mathematics {\bf 1094}, Springer-Verlag (1984) 139--150&

\bibitem[Kis94a]&Kiselman, C.~O.:& Attenuating the singularities of
  plurisubharmonic functions;& Ann.\ Polonici Mathematici {\bf 60} (1994)
  173--197&

\bibitem[Kis94b]&Kiselman, C.~O.:& Plurisubharmonic functions and
  their singularities;& Complex Potential Theory (Eds.\ P.M.~Gauthier \&
  G.~Sabidussi). NATO ASI Series, Series C, Vol.~{\bf 439}
  Kluwer Academic Publishers (1994) 273--323&

\bibitem[OhT87]&T.\ Ohsawa {\rm and} K.\ Takegoshi:& On the extension
  of $L^2$ holomorphic functions;& Math.\ Zeitschrift {\bf 195} (1987)
  197--204&

\bibitem[Puk87]&Pukhlikov, A.V.:& Birational automorphisms of a 
  four-dimensional quintic;& Invent.\ Math.\ {\bf 87} (1987), 303--329&

\bibitem[Puk02]&Pukhlikov, A.V.:& Birational rigid Fano hypersurfaces;& 
  Izv.\ Ross.\ Akad.\ Nauk Ser.\ Mat.\ {\bf 66:}6 (2002), 159--186;
  English translation, Izv.\ Math.\ {\bf 66:} (2002), 1243--1269&

}
\vskip20pt
\noindent
(version of November 23, 2007, printed on \today)
\vfill\eject
\noindent
\section{A. Appendix : a stronger version of Demailly's  estimate
on Monge-Amp\`ere operators}
\medskip
\centerline{\bf Ahmed Zeriahi}
\smallskip
\centerline{Universit\'e Paul Sabatier - Toulouse 3}
\centerline{Laboratoire de Math\'ematiques \'Emile Picard, UMR 5580 du CNRS}
\centerline{118, Route de Narbonne - 31062 Toulouse Cedex 4, France}
\centerline{{\it e-mail\/}: {\tt zeriahi@math.ups-tlse.fr}}
\bigskip

As suggested in Remark 2.10 of J.-P.\ Demailly's paper in the present 
volume, it is possible to weaken the
hypotheses of Theorem~1.1 therein so as to merely assume that the psh function
$\varphi$ on $\Omega$ has zero boundary values on $\partial\Omega$, 
in the sense that the limit of $\varphi(z)$ as $z\in\Omega$ tends to 
any boundary point $z_0\in\partial\Omega$
is zero (see below for an even weaker interpretation). In addition to this,
the integral bound for $e^{-2\varphi}$ can be obtained as a global 
estimate on $\Omega$, and not just on a compact subset~$K\subset\Omega$.
Recall that a complex space is said to be hyperconvex if it possesses
a  bounded (say${}<0$) strictly plurisubharmonic exhaustion function.

\claim (A.1) Theorem|Let $\Omega$ be a bounded hyperconvex domain 
in $\bC^n$ and let $\varphi$ be a 
plurisubharmonic function on $\Omega$ with zero boundary values, such that
$$
\int_\Omega(dd^c\varphi)^n\le M<n^n.
$$
Then there exists a uniform constant $C'(\Omega,M)>0$ independent of $\varphi$
such that
$$
\int_\Omega e^{-2\varphi}d\lambda\le C'(\Omega,M).
$$
\endclaim

\noindent
{\it Proof.} The first step consists of showing that there is a uniform 
estimate
$$
\int_K e^{-2\varphi}d\lambda\le C''(\Omega,K,M)\leqno{\rm(A.2)}
$$
for every compact subset $K\subset\Omega$. Indeed, the compactness 
argument used in the proof of Theorem~1.1 still works in that case, 
thanks to the following observation.

\claim (A.3) Observation|The class $ \cP_{0,M}(\Omega)$ of psh functions 
$\varphi$ on $\Omega$ with zero boun\-dary values and satisfying 
$\int_\Omega (dd^c\varphi)^n\le M$ is a relatively compact subset
of $L^1_\loc(\Omega)$ and its closure $\ol \cP_{0,M}(\Omega)$ consists of
functions sharing the same properties, except that they only have
zero boundary values in the more general sense introduced
by Cegrell {\rm ([Ceg04], see below)}.
\endclaim

This statement is proved in detail in [Zer01]. The argument can be sketched 
as~follows. According to [Ceg04], denote by $\cE_0(\Omega)$ the
set of "test" psh functions, i.e.\ bounded psh functions with zero
boundary values, such that the Monge-Amp\`ere measure has finite mass
on $\Omega$. Then, thanks to $n$ successive integration by parts, one
shows that there exists a constant $c_n>0$ such that if $\varphi$ and
$\psi$ are functions in the class~$\cE_0(\Omega)$, one has 
$$
\int_\Omega(-\varphi)^n (dd^c\psi)^n\le c_n
\Vert\psi\Vert^n_{L^\infty}\int_\Omega(dd^c\varphi)^n.
$$
This estimate is rather standard and was probably stated explicitly for 
the first time by Z.~B{\l}ocki [B{\l}o93].  It is clear by means of 
a standard truncation technique that this estimate is still valid 
when $\varphi\in \cP_{0,M}(\Omega)$ and 
$\psi\in\cE_0(\Omega)$. This proves that $ \cP_{0,M}(\Omega)$ is 
relatively compact in $L^1_\loc(\Omega)$. 

In order to determine the closure $\ol \cP_{0,M}(\Omega)$ of $
\cP_{0,M}(\Omega)$, one can use the class $\cF(\Omega)$ defined by
Cegrell [Ceg04]. By definition, $\cF(\Omega)$ is the class of negative
psh functions $\varphi$ on $\Omega$ such that there exists a non
increasing sequence of test psh functions $(\varphi_j)$ in the class
$\cE_0(\Omega)$ which converges towards $\varphi$ and such that
$\sup_j\int_\Omega (dd^c\varphi_j)^n < +\infty$. Cegrell showed that
the Monge-Amp\`ere operator is still well defined on $\cF(\Omega)$ and is
continuous on non increasing sequences in that space. It is 
then rather easy to show
that the closure of $\cP_{0,M}(\Omega)$ in $L^1_\loc(\Omega)$
coincides with the class of psh functions
$\varphi\in\cF(\Omega)$ such that $\int_\Omega(dd^c\varphi)^n\le
M$. In~fact, if if $(\varphi_j)$ is a sequence of elements of $
\cP_{0,M}(\Omega)$ which converges in $L^1_\loc(\Omega)$ towards
$\varphi$, one knows that $\varphi$ is the upper regularized limit
$\varphi= (\limsup_j \varphi_j)^*$ on~$\Omega$. By putting
$\psi_j:=(\sup_{k\ge j}\varphi_k)^*$, one obtains a non increasing
sequence of functions of $ \cP_{0,M}(\Omega)$ which converges towards
$\varphi$ and since $\varphi_j\le\psi_j\le0$, these functions have 
zero boundary values, and the comparison principle implies that
$\int_\Omega(dd^c\psi_j)^n\le\int_\Omega(dd^c\varphi_j)^n\le M$.  This
proves that $\varphi\in\cF(\Omega)$. The inequality
$\int_\Omega(dd^c\varphi)^n\le M$ also holds true, since $(dd^c\psi_j)^n\to
(dd^c\varphi)^n$ weakly. The estimate (A.2) now follows from the arguments
given by Demailly for Theorem~1.1.

The second step consists in a reduction of Theorem~A.1 to estimate 
(A.2) of the first step, thanks to a subextension theorem with control of the
Monge-Amp\`ere mass. Actually, let $\varphi$ be as above and let
$\tilde\Omega$ be a bounded hyperconvex domain of~$\bC^n$ (e.g.\ a
euclidean ball) such that $\ol\Omega\subset\tilde\Omega$. Then by
[CZ03], there exists $\tilde\varphi\in\cF(\tilde\Omega)$ such that
$\tilde\varphi\le\varphi$ on $\Omega$ and
$\int_{\tilde\Omega}(dd^c\tilde\varphi)^n\le
\int_\Omega(dd^c\varphi)^n\le M$. From this we conclude by (A.2) that
$$
\int_\Omega e^{-2\varphi}d\lambda\le
\int_{\ol\Omega}e^{-2\tilde\varphi}d\lambda\le C''(\tilde\Omega,\ol\Omega,M).
$$
The desired estimate is thus proved with $C'(\Omega,M)=
C''(\tilde\Omega,\ol\Omega,M)$.\qed

\section{References}
\medskip

{\eightpoint

\bibitem[B{\l}o93]&B{\l}ocki, Z.:& Estimates for the complex Monge-Amp\`ere 
operator;& Bull.\ Polish.\ Acad.\ Sci.\ Math.\ {\bf 41} (1993), 151--157&

\bibitem[Ceg04]&Cegrell, U.:& The general definition of the complex 
Monge-Amp\`ere operator;& Ann.\ Inst.\ Fourier {\bf 54} (2004) 159--197&

\bibitem[CZ03]&Cegrell, U., Zeriahi, A.:& Subextension of plurisubharmonic
functions with bounded Monge-Amp\`ere mass;& C.~R.\ Acad.\ Sc.\
Paris, s\'erie I {\bf 336} (2003) 305--308&

\bibitem[Zer01]&Zeriahi, A.:& Volume and Capacity of Sublevel Sets of
a Lelong Class of Plurisubharmonic Functions;& Indiana Univ.\ Math.\ J.\ 
{\bf 50} (2001), 671--703&

}

\end